\documentclass[astronautics,article,submit,pdftex,moreauthors]{Definitions/mdpi}
\firstpage{1} 
\pubvolume{1}
\issuenum{1}
\articlenumber{0}
\pubyear{2026}
\copyrightyear{2025}
\datereceived{ } 
\daterevised{ } 
\dateaccepted{ } 
\datepublished{ } 
\usepackage{mathtools}
\usepackage{bm}
\usepackage[version=4]{mhchem}
\usepackage{siunitx}
\usepackage{algorithm}
\usepackage{algorithmic}
\usepackage{mathrsfs}
\usepackage{bigstrut}
\usepackage{makecell}

\newcommand{\vect}[1]{\boldsymbol{#1}}

\Title{Optimal Control and Neural Porkchop Analysis for Low-Thrust Asteroid Rendezvous Mission}

\Author{Zhong Zhang $^{1,2,}$* \orcidA{}, Niccolò Michelotti $^{1}$, Gonçalo Oliveira Pinho $^{1}$, Yilin Zou $^{2}$ and Francesco Topputo $^{1}$} 

\AuthorNames{Zhong Zhang, Niccolò Michelotti, Gonçalo Oliveira Pinho, Yilin Zou and Francesco Topputo}

\address{%
$^{1}$ \quad Department of Aerospace Science and Technology, Politecnico di Milano, 20156 Milan, Italy

$^{2}$ \quad Department of Aerospace School, Tsinghua University, 100084 Beijing, China}

\corres{Correspondence: zhong.zhang@polimi.it}

\conference{2025 AAS/AIAA Astrodynamics Specialist Conference, Boston} 

\abstract{This paper presents a comparative study of the applicability and accuracy of optimal control methods and neural network-based estimators in the context of porkchop plots for preliminary asteroid rendezvous mission design. The scenario considered involves a deep-space CubeSat equipped with a low-thrust engine, departing from Earth and rendezvousing with a near-Earth asteroid within a three-year launch window. A low-thrust trajectory optimization model is formulated, incorporating variable specific impulse, maximum thrust, and path constraints. The optimal control problem is efficiently solved using Sequential Convex Programming (SCP) combined with a solution continuation strategy.
The neural network framework consists of two models: one predicts the minimum fuel consumption ($\Delta v$), while the other estimates the minimum flight time ($\Delta t$) which is used to assess transfer feasibility.
Case results demonstrate that, in simplified scenarios without path constraints, the neural network approach achieves low relative errors across most of the design space and successfully captures the main structural features of the porkchop plots. In cases where the SCP-based continuation method fails due to the presence of multiple local optima, the neural network still provides smooth and globally consistent predictions, significantly improving the efficiency of early-stage asteroid candidate screening.
However, the deformation of the feasible region caused by path constraints leads to noticeable discrepancies in certain boundary regions, thereby limiting the applicability of the network in detailed mission design phases. Overall, the integration of neural networks with porkchop plot analysis offers an effective decision-making tool for mission designers and planetary scientists, with significant potential for engineering applications.
}

\keyword{asteroid rendezvous; trajectory optimization; low-thrust; neural networks; porkchop plot; optimal control}

\firstpage{1} 
\pubvolume{1}
\issuenum{1}
\articlenumber{0}
\pubyear{2026}
\copyrightyear{2025}
\datereceived{ } 
\daterevised{ } 
\dateaccepted{ } 
\datepublished{ } 

\begin{document}
\nolinenumbers

\section{Introduction}
Near-Earth asteroids (NEAs) offer unique scientific value for understanding the
early evolution of the solar system, assessing extraterrestrial
resources\cite{zhangSustainableAsteroidMining2025b}, and enhancing planetary
defense. They also possess engineering value due to their relative
accessibility. It is estimated that the number of observable NEAs has exceeded
ten thousand, and the number of detectable targets continues to increase
rapidly \cite{granvikDebiasedOrbitAbsolutemagnitude2018}. In recent years,
Japan's Hayabusa-2 mission successfully collected and returned samples from
Ryugu, while NASA's OSIRIS-REx mission obtained samples from Bennu. The DART
mission further demonstrated, for the first time, the feasibility of kinetic
impactor technology for asteroid deflection
\cite{laurettaOSIRISRExSampleReturn2017,chengAIDADARTAsteroid2018,watanabeHayabusa2MissionOverview2017}.
With the maturation of deep-space CubeSats and small satellite platforms,
opportunities for low-cost and low-thrust missions to NEAs are increasing
\cite{ferrariPreliminaryMissionProfile2021a,topputoEnvelopReachableAsteroids2021}.
However, such missions must perform complex trajectory design under strict mass
and power constraints\cite{zhang2025globaloptimalitymultiflybyasteroid},
highlighting the need for more efficient preliminary mission assessment tools.

The ``porkchop'' plot provides an intuitive cost landscape of fuel consumption
and transfer time for interplanetary missions by sweeping departure epochs and
flight durations. However, for low-thrust spacecraft, each grid point requires
solving an optimal control problem that may involve multiple constraints,
variable specific impulse, and other nonlinearities, leading to a computational
burden significantly higher than that of chemical propulsion
\cite{doi:10.2514/2.4231,dellnitzMultiobjectiveApproachDesign2009}. Existing
optimal control methods include direct methods such as pseudospectral methods,
sequential convex programming (SCP), and dynamic programming
\cite{bettsChapter4Optimal2020,zhangGlobalTrajectoryOptimization2024a,malyuta2021convexoptimizationtrajectorygeneration,yamaguchi_trajectory_2025,pavanello_collision_2025},
as well as indirect methods based on Pontryagin’s maximum principle
\cite{jiangPracticalTechniquesLowThrust2012a,bertrandNewSmoothingTechniques2002}.
Even with advanced optimal control methods and continuation techniques
\cite{SCP_formulation}, generating a porkchop plot for a single transfer target
through dense sampling can require tens to hundreds of CPU hours, and is
sensitive to initial guesses and local minima. This remains computationally
prohibitive when applied to a large number of NEA targets.

In recent years, deep learning has been increasingly applied in astrodynamics,
using multilayer perceptrons to directly regress fuel consumption or assess
reachability
\cite{acciariniComputingLowthrustTransfers2024a,zhuFastEvaluationLowThrust2019,liDeepNetworksApproximators2020,10989529,singh_stochastic_2022,kim_probabilistic_2024,10989529,
  MANCINI20233748}. Machine learning models have also been employed for real-time
trajectory guidance
\cite{sanchez-sanchezRealTimeOptimalControl2018,izzoRealTimeGuidanceLowThrust2021,lunghi_energy_2025,smith_model_2025}.
Neural networks, with their universal approximation capabilities
\cite{hornikMultilayerFeedforwardNetworks1989} and microsecond-level inference
efficiency, are well-suited to large-sample scenarios such as porkchop plot
generation. However, existing studies have mostly focused on specific transfer
scenarios, and are often validated under simplified assumptions such as
constant thrust and the absence of path constraints. Moreover, quantitative
comparisons with advanced optimal control methods are often lacking, especially
when addressing realistic engineering constraints such as variable specific
impulse, thermal limits, and power restrictions
\cite{sidhoum2025pontryagin,dambrosioPhysicsInformedPontryaginNeural}.

To this end, this study aims to evaluate the accuracy of neural networks (NNs)
in realistic, constrained low-thrust mission scenarios, and to assess whether
they can serve as a fast surrogate tool for asteroid pre-selection during early
mission design, while also clarifying their limitations. A representative
scenario, porkchop plot generation for rendezvous missions to near-Earth
asteroids, is selected for analysis. The applicability and accuracy of optimal
control methods and neural networks are compared in this context.

The main contributions of this paper are summarized as follows:

\begin{enumerate}
  \item Under realistic scenarios involving variable specific impulse, maximum thrust,
        and path constraints, this paper conducts a quantitative accuracy comparison
        between neural networks and optimal control methods at the porkchop plot level.
  \item Under simple transfer conditions, the neural network maintains a relative error
        of less than 10\% across most of the feasible domain.
  \item The study demonstrates that the porkchop plot is a suitable scenario for
        applying neural networks. It not only provides fast and direct support for
        engineering design, but also allows for visual validation of neural network
        performance across a wide range of transfer scenarios.
\end{enumerate}

The remainder of this paper is organized as follows: Section 2 describes the
mission model and algorithms; Section 3 presents the experimental results;
Section 4 discusses the sources of error and potential improvements; and
Section 5 concludes the paper.

\section{Problem Description and Methods}

This paper investigates the porkchop plot in the context of a practical
asteroid rendezvous mission. The scenario involves a CubeSat launched from
Earth to rendezvous with an asteroid. The rocket provides a maximum hyperbolic
excess velocity $v_{\infty}^{\max}$ of \SI{4}{km/s}, with no constraint on
direction. The departure date ranges from January 1, 2029, to December 31,
2031, with a flight time between 200 and 1200 days and a time step of 10 days.
Each point in the porkchop plot represents a computation aimed at minimizing
the total velocity increment required for the rendezvous. The initial mass of
the satellite is \SI{26.5}{kg}, with \SI{2}{kg} of propellant allocated to the
thruster. The maximum thrust and specific impulse of the low-thrust engine vary
with the spacecraft’s distance from the Sun. Path constraints are also imposed
in this scenario.

This section introduces the optimal control problem corresponding to the
porkchop plot analysis of the asteroid rendezvous mission studied in this
paper. It covers the definition of the low-thrust trajectory optimization
problem, the optimal control method used for its solution, and the neural
network model employed for prediction. First, the mathematical formulation of
the low-thrust trajectory optimization problem is presented, including the
dynamical equations, the modeling of variable specific impulse and maximum
thrust, and the path constraints. Then, the optimal control method used to
solve the problem, namely the Sequential Convex Programming (SCP) algorithm, is
introduced. Finally, the structure and training process of the neural network
model are described, along with the indirect method used for data generation,
and the approach to use the trained model for prediction.

\subsection{Problem Definition}

The problem of interest is the optimization of low-thrust trajectories for
asteroid rendezvous missions. The objective is to minimize the total velocity
increment required for a CubeSat to rendezvous with an asteroid. The equations
of motion governing the spacecraft include gravitational forces and thrust
forces \textcolor{black}{$\vect{\Gamma}$}
\begin{align}
  \dot{\vect{r}} & = \vect{v}                                              \\
  \dot{\vect{v}} & =  - \frac{\mu}{r^3} \vect{r} + \frac{\vect{\Gamma}}{m}
\end{align}
where $\vect{r}$ is the position vector, $\vect{v}$ is the velocity vector, $\mu$ is the gravitational parameter of the Sun, $r = ||\vect{r}||_2$. Fuel consumption obeys
\begin{equation}
  \dot m(t) = -\frac{\|\vect\Gamma(t)\|}{I_{\rm sp}(t)\,g_0}
\end{equation}
where $\vect{\Gamma}(t)$ is the thrust force, $I_{sp}(t)$ is the specific impulse and \textcolor{black}{$g_0$ is the gravitational acceleration at sea level}.

A realistic thruster model is implemented by mapping the maximum thrust
$T_{\max}$ and specific impulse $I_{sp}$ variation over the instantaneous input
power. The available input power $P$ is powered by solar arrays and depends on
the spacecraft distance $r$ from the Sun. \textcolor{black}{Polynomial
  surrogate models of the thruster and power production, derived from
  experimental measurements, are used to capture the complexities of these
  systems while maintaining smooth derivatives.} The coefficients of the
third-order polynomials describing $P(r)$, $T_{\mathrm{max}}(P)$ and
$I_{sp}(P)$ are detailed in Table \ref{tab:Satis_thruster_params}.
\textcolor{black}{ The values correspond to a deep-space CubeSat mission
  comparable to ESA \textit{M-ARGO} mission
  \cite{topputoEnvelopReachableAsteroids2021}. }

\begin{table}[H]
  \caption{Electric propulsion thruster coefficients.}
  \centering
  \begin{tabular}{clll}
    \toprule
    Input power limits & $T_{\mathrm{max}}$ [$\mu$N/W$^i$] coeff. & $I_{sp}$ [s/W$^i$] coeff.    & P [W/AU$^i$] coeff. \\
    \midrule
    \multirowcell{5}{ $P_{\mathrm{in,min}} = 80 \ {\rm W}$                                                             \\ $P_{\mathrm{in,max}} = 130 \ {\rm W}$}
                       & $a_0 = - 1234.3$                         & $b_0 = - 5519.5$             & $c_0 = 1281.4$      \\
                       & $a_1 = 26.498$                           & $b_1 = 225.44$               & $c_1 = - 2518.8$    \\
                       & $a_2 = 0$                                & $b_2 = - 1.8554$             & $c_2 = 1828.5$      \\
                       & $a_3 = 0$                                & $b_3 = 5.0836\times 10^{-3}$ & $c_3 = -475.80$     \\
    \bottomrule
  \end{tabular}
  \label{tab:Satis_thruster_params}
\end{table}

\begin{figure}[H]
  \centering
  \includegraphics[width=0.6\linewidth]{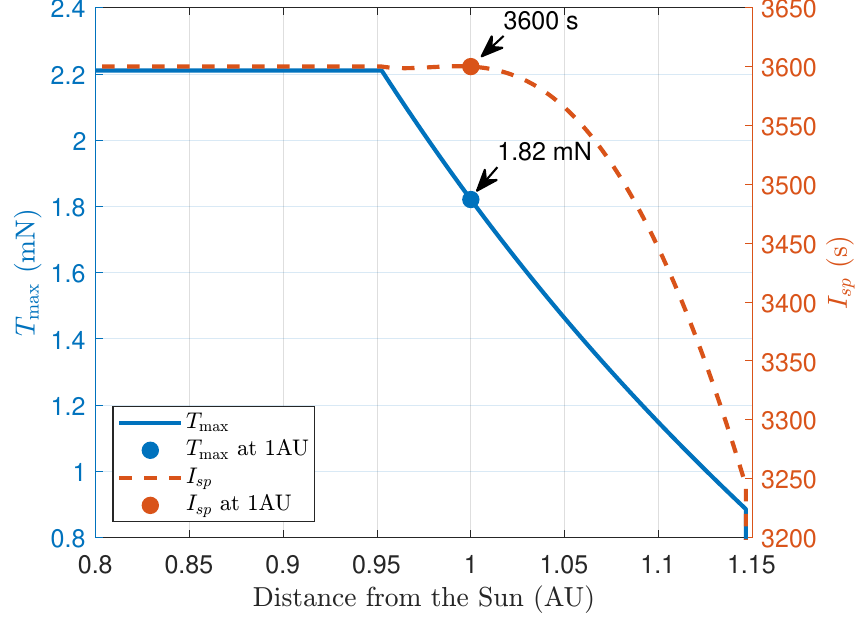}
  \caption{Variable $T_{\rm max}$ and $I_{\rm sp}$ curves with respect to the distance from the Sun.}
  \label{fig:variable_Tmax_Isp}
\end{figure}

Maximum thrust and specific impulse can also be expressed directly as a
function of the distance from the Sun, as shown in Figure
\ref{fig:variable_Tmax_Isp}.

Variations of \(T_{\max}\) and \(I_{\rm sp}\) with heliocentric distance (in
astronomical units, AU) are modeled by fourth order polynomials between 0.9525
AU to 1.1467 AU
\begin{equation}
  \textcolor{black}{
    T_{\max}(x)}=
  \begin{cases}
    2.2\times10^{-3}\ \mathrm{N},           & x\le0.9525,      \\[6pt]
    \displaystyle\sum_{i=0}^4 c_{T,i}\,x^i, & 0.9525<x<1.1467, \\[6pt]
    0,                                      & x\ge1.1467,
  \end{cases}
  \quad
  I_{\rm sp}(x)=
  \begin{cases}
    3600\ \mathrm{s},                       & x\le0.9525,      \\[6pt]
    \displaystyle\sum_{i=0}^4 c_{I,i}\,x^i, & 0.9525<x<1.1467, \\[6pt]
    0,                                      & x\ge1.1467,
  \end{cases}
  \label{Eq.isp}
\end{equation}
where the thrust coefficients and the \(I_{\rm sp}\) coefficients are summarized in Table \ref{tab:coeffs}.

\begin{table}[H]
  \caption{Fourth‐order polynomial coefficients for thrust $T_{\max}(x)$ and specific impulse $I_{\rm sp}(x)$.}
  \label{tab:coeffs}
  \centering
  \begin{tabular}{c c c}
    \toprule
    $i$ & $c_{T,i}$       & $c_{I,i}$                   \\
        & [N\,/\,AU$^i$]  & [s\,/\,AU$^i$]              \\
    \midrule
    0   & $0.0327202372$  & $3.1951149396 \times10^{5}$ \\
    1   & $-0.0667431624$ & $-1.1891014516\times10^{6}$ \\
    2   & $0.0484515930$  & $1.6645687780 \times10^{6}$ \\
    3   & $-0.0126077484$ & $-1.0253620793\times10^{6}$ \\
    4   & $0$             & $2.3398323687 \times10^{5}$ \\
    \bottomrule
  \end{tabular}
\end{table}

Due to thermal and power constraints, the CubeSat trajectory is required to
stay above 0.8 AU, and the engine must be shut down beyond 1.15 AU throughout
the mission, representing a path constraint.

\subsection{Optimal Control Method}

To efficiently analyze the reachability of multiple asteroids, a fast and
reliable method is needed to compute thousands of fuel-optimal low-thrust
trajectories per target. This work employs a direct Sequential Convex
Programming (SCP) algorithm, solving a convex subproblem for each combination
of departure date and time of flight. The solution is then used as an initial
guess for neighbouring problems across the time grid, ensuring robust
convergence and smooth porkchop plots. The nonconvex low-thrust trajectory
optimization problem is solved by considering a sequence of convex subproblems
whose solution converges, under defined hypotheses, to the solution of the
original one~\citep{malyuta2021convexoptimizationtrajectorygeneration}.

\textcolor{black}{
At each SCP iteration, the quantities $\tilde{\vect{x}}(t)$ and $\tilde{\vect{u}}(t)$ denote the reference state and control vectors obtained from the solution of the previous iteration. These reference trajectories are used to build a local convex approximation of the original nonconvex problem through first-order linearization. Exploiting a correct change of variables and convexification techniques, at each iteration of the SCP, the following problem is solved:
\begin{subequations} \label{eq:opt_problem}
  \begin{align}
    \min_{\vect{u}(t)} \quad & -\kappa(t_f) + \lambda_1 \int_{t_i}^{t_f} \max(0, \eta(t))\, dt + \lambda_2 \int_{t_i}^{t_f} \| \vect{\gamma}(t) \|_1\, dt \label{eq:opt_problem_objective} \\
    \text{s.t.:} \quad
                             & \dot{\vect{x}}(t) =
    \vect{f}(\tilde{\vect{x}}(t), \tilde{\vect{u}}(t))
    + \vect{A}(\tilde{\vect{x}}(t)) \big(\vect{x}(t) - \tilde{\vect{x}}(t)\big)
    + \vect{B} \big( \vect{u}(t) \big)
    + \vect{\gamma}(t) \label{eq:opt_problem_dyn}                                                                                                                                          \\
                             & \tau_x^2 + \tau_y^2 + \tau_z^2 \leq \tau \label{eq:opt_problem_v_boun d}                                                                                    \\
                             & \tau(t) \leq T_{\max} e^{-\tilde{\kappa}} [1 - \kappa +\tilde{\kappa}] + \eta(t) \label{eq:opt_problem_gamma_bound}                                         \\
                             & \| \vect{x}(t) - \tilde{\vect{x}}(t) \|_1 \leq R \label{eq:opt_problem_state_dev}                                                                           \\
                             & \vect{r}(t_i) = \vect{r}_i, \quad \vect{v}(t_i) = \vect{v}_i, \quad \kappa(t_i) = \kappa_i \label{eq:opt_problem_init}                                      \\
                             & \vect{r}(t_f) = \vect{r}_f, \quad \vect{v}(t_f) = \vect{v}_f \label{eq:opt_problem_final}                                                                   \\
                             & \vect{x}_l \leq \vect{x}(t) \leq \vect{x}_u, \quad \vect{u}_l \leq \vect{u}(t) \leq \vect{u}_u \label{eq:opt_problem_bounds}
  \end{align}
\end{subequations}
where $\vect{x} = [\vect{r}, \vect{v}, \kappa]$ and $\vect{u} = [\tau_x, \tau_y, \tau_z, \tau]$ are the state and control variables, respectively. A change of variables has been exploited by defining $\kappa = \ln{m}$ and $ \tau = \vect{\Gamma}/m$ ($\tau_x, \tau_y, \text{and } \tau_z$ defined accordingly). The slack variables $\eta(t)$ and $\vect{\gamma}(t)$ are introduced to avoid artificial infeasibility of the convex subproblem: $\eta(t)$ relaxes the thrust bound in Eq.~\eqref{eq:opt_problem_gamma_bound}, while $\vect{\gamma}(t)$ relaxes the linearized dynamics in Eq.~\eqref{eq:opt_problem_dyn}. The penalty parameters $\lambda_1$ and $\lambda_2$ are selected sufficiently large to enforce $\eta(t)\le 0$ and $\vect{\gamma}(t)\approx \vect{0}$ at convergence.
The thrust bound uses the distance-dependent maximum thrust $T_{\max}(r)$
defined in Eq.~\eqref{Eq.isp}, and $\eta(t)$ is introduced to avoid artificial
infeasibility in the convex subproblem. The trust-region radius $R$ limits the
deviation from the reference trajectory and preserves the validity of the local
linearization.}

\textcolor{black}{
  In Eq.~\eqref{eq:opt_problem_dyn}, the vector of dynamics $\vect{f}$ and matrices $\vect{A}$ and $\vect{B}$ are defined as
  \begin{equation}
    \mathbf{f}\bigl(\mathbf{x}(t), \mathbf{u}(t)\bigr)
    =
    \begin{bmatrix}
      \mathbf{v}(t)                      \\
      -\mu \, \dfrac{\mathbf{r}(t)}{r^3} \\
      0
    \end{bmatrix}
    +
    \underbrace{
      \begin{bmatrix}
        \mathbf{0}_{3 \times 3} \;\; \mathbf{0}_{3 \times 1} \\
        \mathbf{I}_{3 \times 3} \;\; \mathbf{0}_{3 \times 1} \\
        \mathbf{0}_{1 \times 3} \;\; -\dfrac{1}{g_0 I_{sp}}
      \end{bmatrix}
    }_{=:{\mathbf{B}}}
    \begin{bmatrix}
      \tau_x(t) \\
      \tau_y(t) \\
      \tau_z(t) \\
      \tau(t)
    \end{bmatrix},
  \end{equation}
  \begin{equation}
    \vect{A}(\tilde{\vect{x}}(t)) \coloneqq \left. \frac{\partial \vect{f}}{\partial \vect{x}} \right|_{\tilde{\vect{x}}(t),\,\tilde{\vect{u}}(t)},
  \end{equation}
  where $r = ||\vect{r}(t)||_2$.
}
In Eq.~\eqref{eq:opt_problem_init}, $\vect{v}_i$ is given by the sum between the heliocentric velocity of the Earth at $t_i$ and the initial velocity provided by the launcher $\mathbf{v}_\infty$, which translates as
\begin{equation}
  \left\{
  \begin{array}{l}
    \vect v_i = \vect v_E(t_i)  + \vect v_\infty \\
    || \vect v_\infty ||_2  \le v_\infty^\mathrm{max}
  \end{array}
  \right.
  \label{eq:ini_velocity}
\end{equation}
where $v_\infty^\mathrm{max}$ is equal to \SI{4}{km/s}. Note that the ecliptic angles of the vector $v_\infty^\mathrm{max}$, namely the right ascension and declination of the launcher velocity, are not constrained.

The problem in Eqs. \eqref{eq:opt_problem} is a Second-Order Cone Program
(SOCP) and can be solved by efficient convex solvers. In this work, the
Embedded COnic Solver (ECOS) \cite{Domahidi2013ecos} has been used.
Additionally, an arbitrary-order Legende-Gauss-Lobatto method based on Hermite
interpolation is used \cite{HS_discretization}. Despite the robustness of the
SCP algorithm, the solution of a trajectory optimization problem depends on the
initial guess provided. Therefore, to obtain a homogenous porkchop plot a
\textit{continuation method} has been used. It consists of exploiting the
solutions of already-solved problems and using them as initial guesses for
their neighbours. The continuation scheme adopted follows the one implemented
in Morelli et al.~(2024) \cite{MORELLI20244241}. In particular, the solution at
earliest departure date and maximum time of flight (which corresponds to the
top-left corner of the porkchop plot) is solved by trying two different initial
guesses: a simple shape-based method approximated by a cubic polynomial
function \cite{Shape-Based_IG} and a more sophisticated sub-optimal trajectory
resulting from an adapted version of the Q-law guidance algorithm \cite{q_law}.
Then, this solution is used as the initial guess of lower time of flights and
later departure dates. Each time a new problem is solved, its solution becomes
the initial guess for neighbour ones.

The idea behind continuation is that, if the search space of a porkchop plot is
discretized with a fine enough time grid, the nearby problems are characterized
by fairly similar initial and final conditions, and therefore the solution of
one represents a good initial guess for the other. Still, artificial
discontinuities in the solutions space can be introduced by mission
assumptions, such as the spacecraft thermal limits or thruster model
constraints. In particular, the spacecraft thermal limits introduce a
hard-constraint on the minimum distance from the Sun of 0.8 AU, below which the
spacecraft cannot go. In addition, the cut-off threshold of the input power for
the thruster forces the shutdown of the engine at distances above approximately
1.15 AU. \textcolor{black}{There may also be bifurcations or island artifacts
  in the solution space that cannot be effectively handled by a continuation
  approach.} Such situations can become relevant, increasing the complexity of
the problem and leading to local optimal convergence.

To tackle these problems when encountered, the gradual introduction of the
constraints is exploited through a two-step relaxation strategy. Namely, the
porkchop plot is computed first without some constraints to ensure a smooth
solution space. Afterward, each solution from the relaxed problem is used as
the initial guess of the refined porkchop plot, which considers the original
constrained problem. Although this approach is undoubtedly time-consuming, it
improves the convergence and robustness when computing porkchop plots. However,
in a limited number of instances, some discontinuities could remain visible in
the porkchop plots. This is not necessarily due to a convergence problem of the
solver and bifurcations in the solution space led by the continuation approach;
rather, it could be related to the actual infeasibility of solutions in certain
combinations of time of flight and departure date. Such combinations indeed
determine the relative configuration between the departure and arrival points,
which may lead to infeasible trajectories once the thruster model and distance
constraints are injected into the problem.

\subsection{Neural Network Approximator}

The neural network is trained on data generated by an indirect method solver
with no path constraints and with constant specific impulse $I_{\rm sp}$ and
maximum thrust $T_{\max}$. The neural networks predict the optimal transfer
time and the minimum fuel consumption. In this subsection, we first describe
how the indirect method is used to generate the training data, then we present
the architecture of the neural network, and finally we explain how the trained
network is used to produce pork‑chop diagrams.

\subsubsection{Indirect Method for Data Generation}

Define the normalized thrust control \(\vect u(t)\), with \(\|\vect u\|\le1\),
and thrust vector defined as
\[
  \vect\Gamma = T_{\max}(r)\,\vect u.
\]
Instead of using Cartesian coordinates \((\boldsymbol{r}, \boldsymbol{v})\) or
classical orbital elements \((a, e, i, \Omega, \omega, f)\), the state vector
\(\boldsymbol{S}\) is represented using modified equinoctial elements (MEE)
\cite{walkerSetModifiedEquinoctial1985,junkinsExplorationAlternativeState2019},
which are non-singular and well-suited for solving the majority of low-thrust
trajectory optimization problems.
\begin{equation}
  \begin{array}{ccc}
    p = a (1 - e^2),            & f = e \cos (\omega + \Omega), & g = e \sin (\omega + \Omega),   \\
    h = \tan (i/2) \cos \Omega, & k = \tan (i/2) \sin \Omega,   & L = \omega + \Omega + \upsilon,
  \end{array}
\end{equation}
where \(a\), \(e\), \(i\), \(\Omega\), \(\omega\), and \(\upsilon\) represent the semi-major axis, eccentricity, inclination, right ascension of the ascending node, argument of perigee, and true anomaly, respectively. Defining the spatial state vector as \(\bm{s} = [p,\,f,\,g,\,h,\,k,\,L]^{\mathrm{T}}\), the dynamics are then expressed as
\begin{equation}
  \label{eq:indirect_dynamics}
  \dot{\bm{s}} = \bm{D}(\bm{s}) + \frac{T_{\max}(r)}{m}\,\bm{M}(\bm{s})\,\bm{\alpha}\,u.
\end{equation}
\textcolor{black}{
  where $u \in [0,1]$ denotes the normalized thrust magnitude, and $\boldsymbol{\alpha}$ is the unit thrust direction vector, satisfying $\lVert \boldsymbol{\alpha} \rVert = 1$.}
Detailed expressions of \(\bm{M}\) and \(\bm{D}\) are provided in \cite{izzoRealTimeGuidanceLowThrust2021,junkinsExplorationAlternativeState2019,gaoLowThrustInterplanetaryOrbit2004}. Two optimization objectives are considered:
\begin{itemize}
  \item {Fuel-Optimal:} Minimize
        \begin{equation}
          \label{eq:indirect_log_homotopy}
          J_{\rm fuel}  =
          \lambda_0 \int_{t_0}^{t_f} \frac{T_{\max }}{I_{\mathrm{sp}} g_0}  \{u-\varepsilon \ln [u(1-u)]\} \mathrm{d} t
        \end{equation}
        subject to the dynamics \(\dot{\vect s}, \dot m\) and boundary conditions
        \(\vect s(t_0)=\vect s_0\), \(\vect s(t_f)=\vect s_f\), \(m(t_0)=m_0\). \textcolor{black}{$\varepsilon > 0$ is a small smoothing parameter introducing a logarithmic barrier to prevent singularities at $u=0$ and $u=1$, thereby improving numerical convergence\cite{jiangPracticalTechniquesLowThrust2012a} .}
  \item {Time-Optimal:} Minimize
        \begin{equation}
          \label{eq:indirect_log_homotopy_time}
          J_{\rm time} =
          \lambda_0 \int_{t_0}^{t_f} 1 \mathrm{d} t
        \end{equation}
        subject to the same dynamics and endpoint constraints, with free final time \(t_f\).
\end{itemize}

The Hamiltonian for the fuel problem is given by
\begin{adjustwidth}{-\extralength}{0cm}
  \begin{equation}
    H_{\rm fuel}
    = \bm{\lambda}_s^{\mathrm{T}} \dot{\bm{s}}
    + \lambda_m \dot{m}
    + \lambda_0 L_{\rm fuel}(\bm{s},\bm{u})
    = \frac{T_{\max}}{m} \bm{\lambda}_s^{\mathrm{T}} \bm{M}\,\bm{\alpha}\,u
    + \bm{\lambda}_s^{\mathrm{T}} \bm{D}
    - \lambda_m \frac{T_{\max}}{I_{\rm sp}\,g_0}\,\bm{\alpha}\,u
    + \lambda_0 \frac{T_{\max}}{I_{\rm sp}\,g_0}\bigl[u - \varepsilon \ln\bigl(u(1-u)\bigr)\bigr].
  \end{equation}
\end{adjustwidth}
Costate dynamics are given by
\[
  -\dot{\vect \lambda} = \frac{\partial H}{\partial(\vect s,m)}.
\]
Pontryagin’s minimum principle yields
\[
  \vect u^* = -\frac{\vect\Phi}{\|\vect\Phi\|},
  \quad
  \vect\Phi \equiv \frac{T_{\max}(r)}{m}\,\vect\lambda_v
  - \lambda_m\,\frac{T_{\max}(r)}{I_{\rm sp}(r)\,g_0}\,\frac{\vect\lambda_v}{\|\vect\lambda_v\|}.
\]
The final-time transversality condition for free \(m(t_f)\) gives
\(\lambda_m(t_f)=0\). Normalization of the costate norm \(\|\vect\lambda(t_0)\|
= 1\) fixes the scaling. For time-optimal control, the running cost reduces to
\(\lambda_0\), the control is always full throttle \((\|\vect u\|=1)\), and the
free-final-time condition \(H(t_f)=0\) determines \(t_f\). The time-optimal
Hamiltonian becomes
\begin{equation}
  H_{\rm time}
  = \frac{T_{\max}}{m} \bm{\lambda}_s^{\mathrm{T}} \bm{M}\,\bm{\alpha}\,u
  + \bm{\lambda}_s^{\mathrm{T}} \bm{D}
  - \lambda_m \frac{T_{\max}}{I_{\rm sp}\,g_0}\,\bm{\alpha}\,u
  + \lambda_0.
\end{equation}

The resulting two-point boundary value problem is solved using a single
shooting method, yielding the optimal state trajectories \(\vect r(t)\),
\(\vect v(t)\), \(m(t)\), and thrust history \(\vect\Gamma(t)\).
\begin{itemize}
  \item {Fuel-Optimal Shooting Function:}
        Let the state be \(\vect x=[\vect s;\,m]\) with costate \(\vect\lambda=[\lambda_{\vect r};\,\lambda_{\vect v};\,\lambda_m]\) and multiplier \(\lambda_0\). The unknowns are the initial costates \(\vect\lambda(t_0)\). The shooting function is
        \begin{equation}
          \Phi_{\rm fuel}\bigl[\vect\lambda(t_0)\bigr]
          = \begin{bmatrix}
            \vect s(t_f)-\vect s_f \\[4pt]
            \lambda_m(t_f)         \\[4pt]
            \|\vect\lambda(t_0)\|-1
          \end{bmatrix}
          = \vect 0.
        \end{equation}

  \item {Time-Optimal Shooting Function:}
        Let \(t_f\) be free and the terminal state depend on \(t_f\). The unknowns are \(\vect\lambda(t_0)\) and \(t_f\). The boundary value function is
        \begin{equation}
          \Phi_{\rm time}\bigl[\vect\lambda(t_0);t_f\bigr]
          = \begin{bmatrix}
            \vect s(t_f)-\vect s_f(t_f)                                     \\[4pt]
            \lambda_m(t_f)                                                  \\[4pt]
            H_{\rm time}(t_f)-\vect\lambda_L(t_f)\,\frac{\sqrt{\mu p}}{r^2} \\[4pt]
            \|\vect\lambda(t_0)\|-1
          \end{bmatrix}
          = \vect 0.
        \end{equation}
\end{itemize}

\subsubsection{Neural Network Model}
This paper focuses on the application and comparison of neural networks and
optimal control methods in the analysis of asteroid rendezvous missions. The
training procedure and performance evaluation of the neural network are not the
primary concerns of this study; therefore, a pre-trained model is adopted
directly. Relevant work can be found in \cite{zhang_neural_2025}.

The inputs to the neural network interface include the departure and arrival
positions and velocities, the transfer time, the maximum thrust, the specific
impulse, and the initial mass. Two separate models are used for the outputs:
one predicts the solution to the $\Delta v$ fuel-optimal problem, and the other
predicts the $\Delta t$ solution to the time-optimal problem. The time-optimal
$\Delta t$ model is used to assess reachability—if the given transfer time is
shorter than the predicted minimum transfer time, the transfer is considered
infeasible. The final deep neural network employed for both the \(\Delta v\)
and \(\Delta t\) models consists of \(n_{\rm layer}=9\) fully connected hidden
layers, each containing \(n_{\rm neuron}=128\) neurons. The ReLU activation
function is used throughout all hidden layers to introduce nonlinearity, while
the output layer applies a linear activation to directly predict the continuous
transfer cost. Training was performed on a 100-million-sample dataset, which
requires approximately three days to produce each 1-million-sample block,
culminating in a total generation time of three days for the full dataset.
Model training, executed on the same computational platform, also completes in
about three days. On the test sets, the \(\Delta v\) model achieves a mean
absolute error (MAE) of \SI{3.38}{m/s} and a mean relative error (MRE) of
0.78\%, significantly outperforming prior methods
\cite{zhuFastEvaluationLowThrust2019} \cite{liDeepNetworksApproximators2020}
\cite{acciariniComputingLowthrustTransfers2024a}. For the \(\Delta t\) model,
the MAE on minimum transfer time is 2.56 days, with an MRE of 0.63\%.
Additionally, the neural network model leverages dynamic transformations, and
the generated dataset spans a maximum acceleration range of
\(2.5\times10^{-6}\) to \(1.2\times10^{-2}\)\,m/s\(^2\) and a specific impulse
range of 700 to 9000\,s, enabling the network’s strong generalization across
arbitrary orbital parameters (\(a\), \(e\), \(i\)) and propulsion
configurations (\(a_s\), \(I_{\rm sp}\)) to support almost all current
low-thrust missions without requiring dataset retraining.

\subsubsection{Multi-revolution Trajectory Optimization}


The neural network model considered here is applicable only to
single-revolution transfers and cannot be directly extended to trajectories
involving multiple revolutions. To address this limitation, a hybrid
optimization framework is introduced that combines global search, local
refinement, and feasibility checking based on the \(\Delta t\) model. First,
the multi-revolution trajectory is uniformly partitioned into several time
segments, each spanning less than one orbital revolution. The optimization
variables include the orbital states at the intermediate waypoints between
adjacent segments. In addition, three further optimization variables correspond
to the components of the launch velocity provided by the launch vehicle. The
optimization process consists of two stages: a global search using particle
swarm optimization, followed by local refinement using a gradient-based
optimization algorithm. The objective function is defined as the total velocity
increment over all segments. Moreover, each segment incorporates a feasibility
constraint evaluated by an additional neural network model and enforced through
a penalty function.

\begin{figure}[H]
  \centering
  \includegraphics[width=0.4\linewidth]{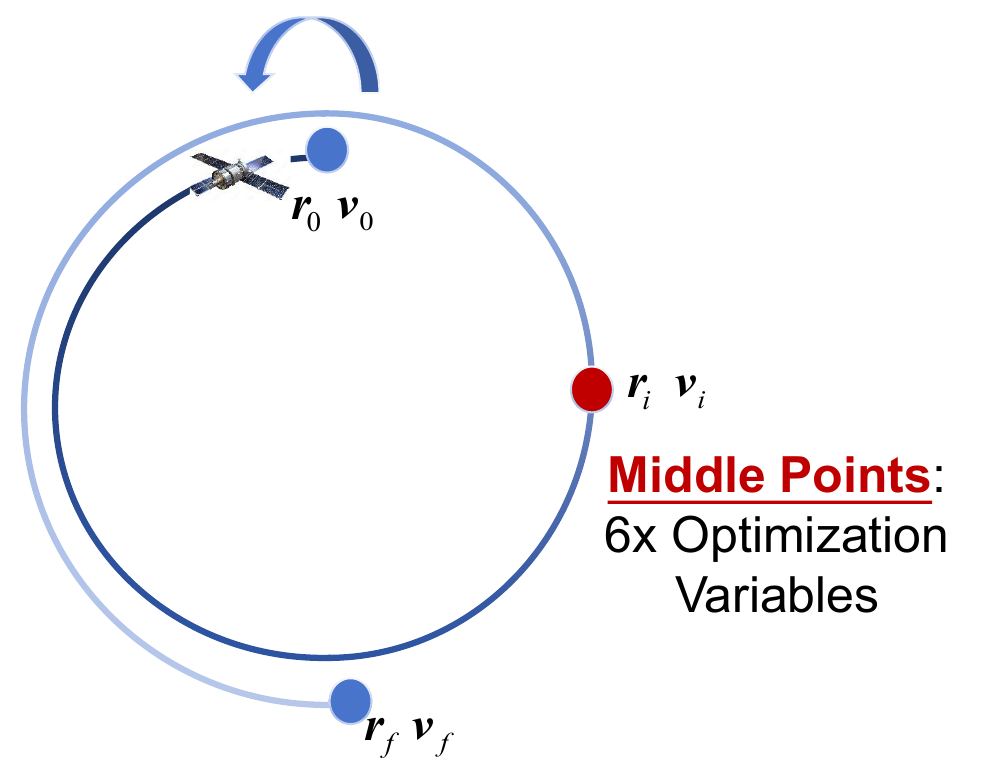}
  \caption{Connecting trajectories for computing multi-revolution transfers.}
  \label{fig:connecting}
\end{figure}

The detailed procedure is as follows: First, the full multi‐revolution transfer
is partitioned into \(N\) fixed time segments
\(\{[t_{k},t_{k+1}]\}_{k=0}^{N-1}\), each spanning less than one orbital
period. The optimization variables are:
\[
  \vect s_k = [\,\vect r(t_k),\,\vect v(t_k)\,]_{k=1}^{N-1},
  \quad
  \Delta \bm{v}_{\rm launch} = [\Delta v_x,\,\Delta v_y,\,\Delta v_z]
\]
where \(\vect r(t_k),\,\vect v(t_k)\) denotes the spatial state at the \(k\)-th
intermediate waypoint, and \(\Delta \bm{v}_{\rm launch}\) represents the three
components of the initial injection velocity supplied by the launch vehicle.
The specific impulse $I_{\rm sp}$ and maximum thrust $T_{\max}$ can either be
fixed or vary at each intermediate waypoint according to Eq.~\ref{Eq.isp}, the
latter case being referred to as "piece-wise variable". The total velocity
increment objective is
\[
  J = \sum_{k=0}^{N-1} \Delta v_k^\ast
\]
where \(\Delta v_k^\ast\) is the optimal velocity increment for segment \(k\),
as predicted by our single‐revolution neural \(\Delta v\) model. To ensure
dynamic consistency and mission feasibility, each segment’s end state must
satisfy both boundary matching and feasibility constraints. The latter are
imposed via a penalty function
\[
  P_{\rm feas}(\bm{s}_k, \bm{s}_{k+1})
  = \rho \,\max\bigl\{0,\;
  (t_{k+1} - t_k) - \mathcal{N}_{\Delta t}(\bm{s}_k,\bm{s}_{k+1}) \bigr\}^2,
\]
where \(\mathcal{N}_{\Delta t}\) is the output of the time-optimal \(\Delta t\)
model, with \(\rho\) a large penalty weight. The overall optimization proceeds
in two stages:
\begin{enumerate}
  \item {Global Search:} A particle swarm optimization (PSO) algorithm explores the high‐dimensional optimization space \(\{\bm{s}_k,\Delta \bm{v}_{\rm launch}\}\), minimizing \(J + \sum_k P_{\rm feas}\).
  \item {Local Refinement:} The best candidate from PSO is passed to a gradient‐based local solver \cite{Johnson2010} to fine‐tune the intermediate states and launch velocity, further reducing the objective while strictly enforcing segment continuity and feasibility.
\end{enumerate}

\section{Results}

\begin{figure}[H]
  \centering
  \includegraphics[width=0.7\linewidth]{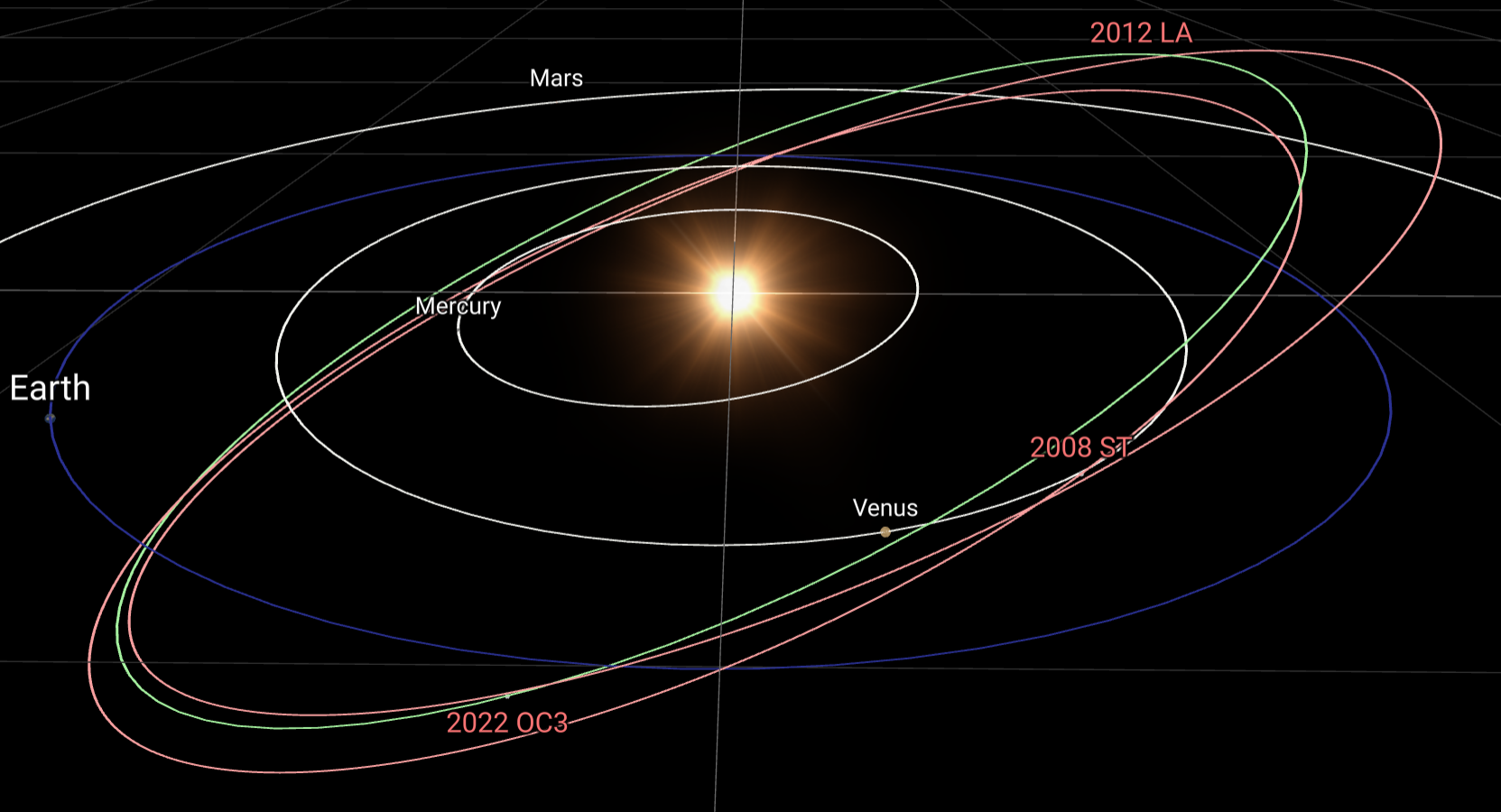}
  \caption{Orbits of selected asteroids. Credits: ESA NEO Toolkit}
  \label{fig:orbital_elements}
\end{figure}

Among the numerous near-Earth asteroid analyzed, three representative porkchop
plots, corresponding to asteroids 2012 LA, 2008 ST, and 2022 OC3, are presented
for demonstration. Their orbital elements are shown in Table
\ref{tab:orbital_elements} and Figure \ref{fig:orbital_elements}.

\begin{table}[H]
  \caption{Classical orbital elements of selected asteroids.}
  \label{tab:orbital_elements}
  \centering
  \begin{tabular}{@{}c l c c c c c c@{}}
    \toprule
    Epoch & Asteroid & $a$ (km)              & $e$       & $i$ (°) & $\Omega$ (°) & $\omega$ (°) & $\nu$ (°) \\
    \midrule
    \multirow{3}{*}{ \makecell{2029 Jan 01                                                                   \\12:00:00 UTC}}
          & 2012 LA  & $1.55604\times10^{8}$ & 0.0216877 & 21.7262 & 353.0774     & 320.1778     & 257.0437  \\
          & 2008 ST  & $1.43986\times10^{8}$ & 0.1265623 & 21.5865 & 359.3513     & 120.5766     & 157.4352  \\
          & 2022 OC3 & $1.52899\times10^{8}$ & 0.1286078 & 22.9596 & 0.5957       & 49.6797
          & 113.2483                                                                                         \\
    \bottomrule
  \end{tabular}
\end{table}

\begin{figure}[H]
  \centering
  \includegraphics[width=0.8\linewidth]{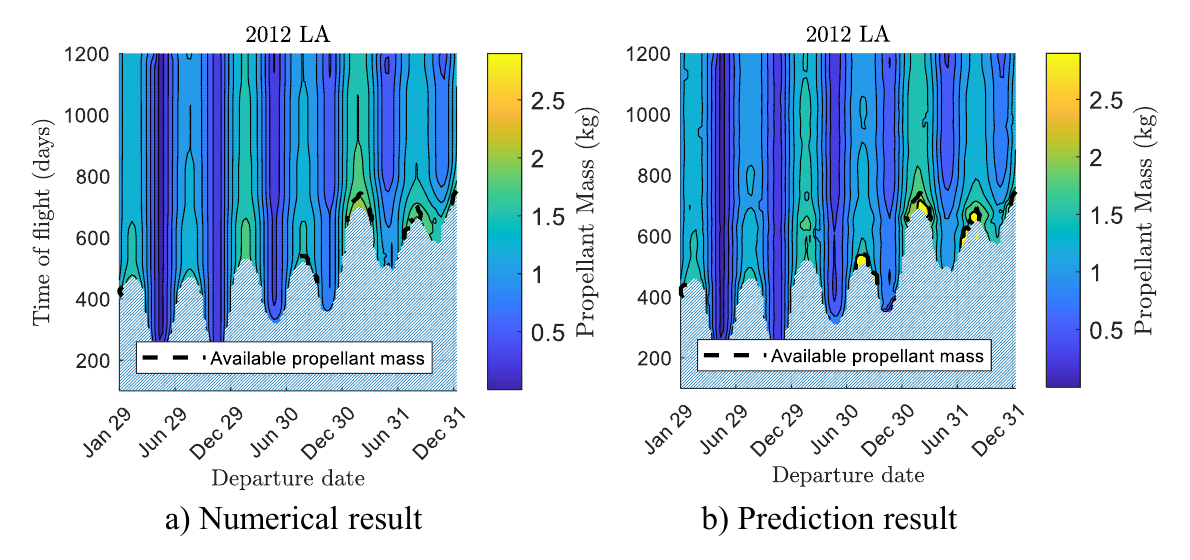}
  \caption{Porkchop plot for the 2012 LA Asteroid.}
  \label{fig:porkchop_2012LA}
\end{figure}

\begin{figure}[H]
  \centering
  \includegraphics[width=0.8\linewidth]{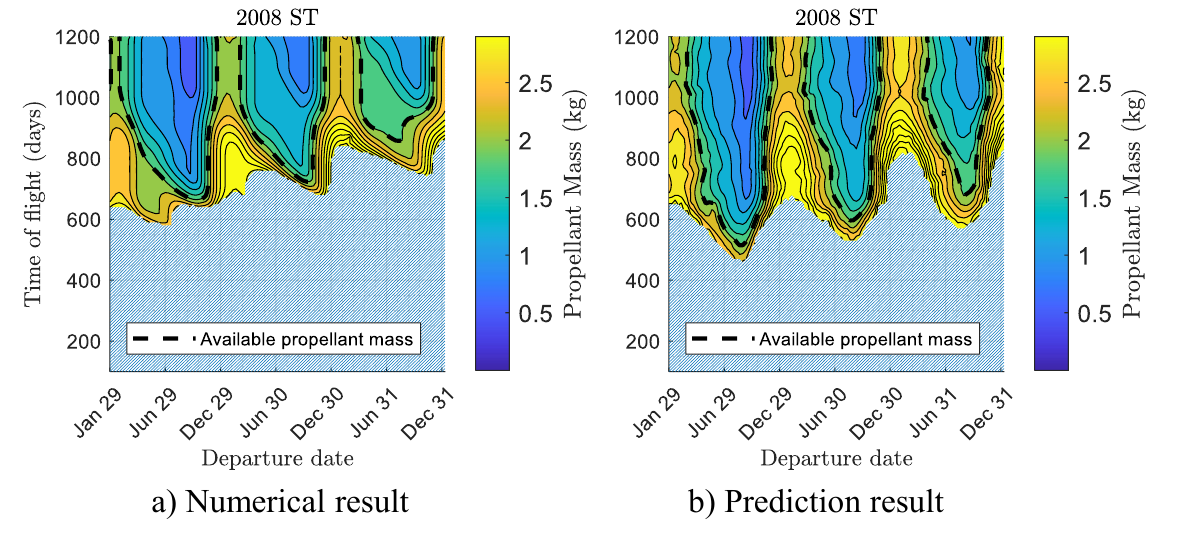}
  \caption{Porkchop plot for the 2008 ST Asteroid.}
  \label{fig:porkchop_2008ST}
\end{figure}

\begin{figure}[H]
  \centering
  \includegraphics[width=0.8\linewidth]{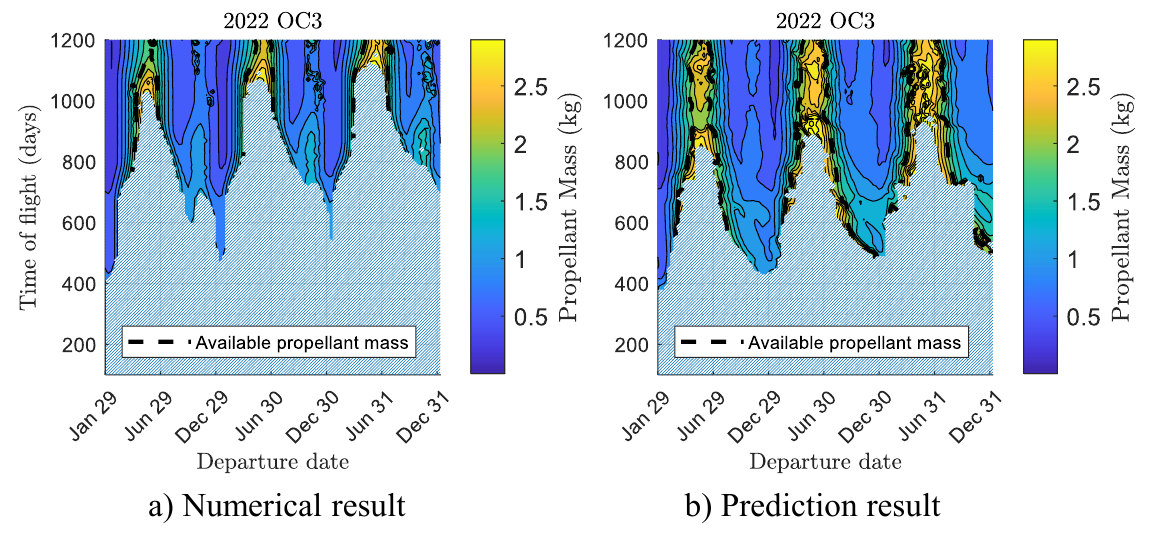}
  \caption{Porkchop plot for the 2022 OC3 Asteroid.}
  \label{fig:porkchop_2022OC3}
\end{figure}

The porkchop plots for these three asteroids are shown in Figures
\ref{fig:porkchop_2012LA}, \ref{fig:porkchop_2008ST}, and
\ref{fig:porkchop_2022OC3}. In each plot, the horizontal axis represents the
launch date, the vertical axis represents the time of flight, and the color
scale indicates fuel consumption. For comparison, the direct method results are
compared with those from the neural network model. The direct method (numerical
results) accounts for variable specific impulse, maximum thrust, and includes
path constraints, whereas the neural network model (prediction results) is
based on an indirect method using piece-wise variable specific impulse and
maximum thrust. For asteroid 2012 LA, large regions of low fuel consumption
appear in the porkchop plot, and the neural network model results closely match
the direct method. The porkchop plot for asteroid 2008 ST exhibits a more
complex structure, with multiple high fuel consumption regions, and significant
discrepancies between the neural network model results and the direct method in
some areas. The porkchop plot for asteroid 2022 OC3 exhibits even more
pronounced discrepancies, as the optimal control direct method results display
numerous discontinuities and gaps. The characteristics of these distinct
porkchop plots will be discussed in detail in the following section.

\section{Analysis and Discussion}

\begin{figure}[H]
  \centering
  \includegraphics[width=1.0\linewidth]{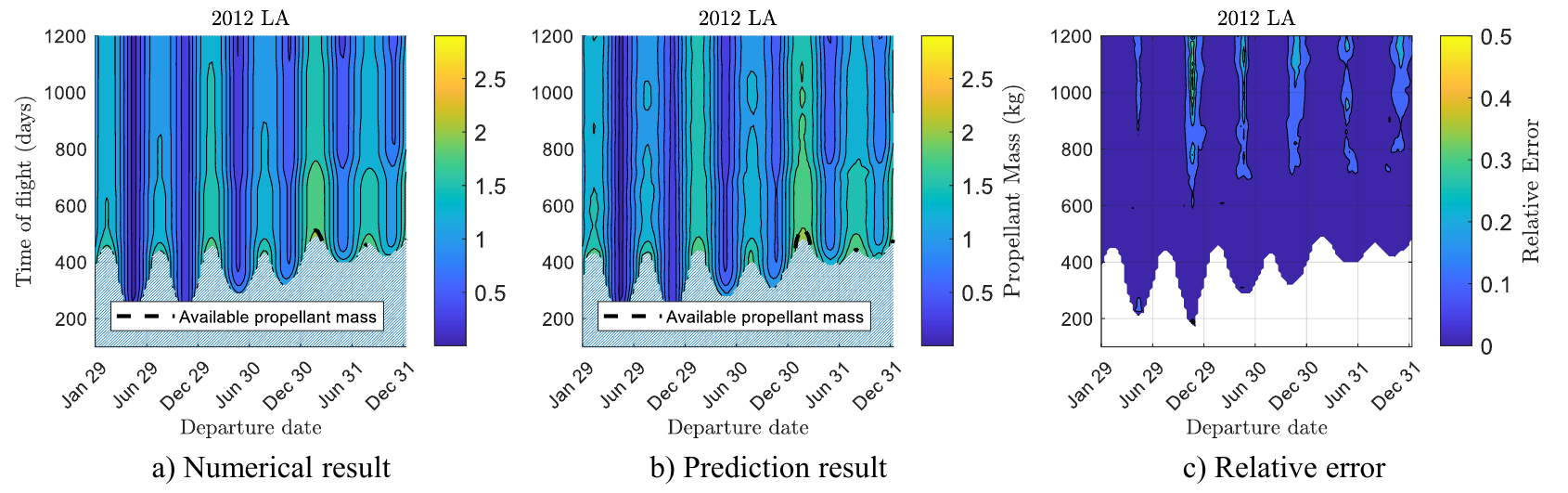}
  \caption{Porkchop plot for the Asteroid 2012 LA without variable $I_{\rm sp}$, $T_{\max}$, and path constraints.}
  \label{fig:2012LA_noconstraint}
\end{figure}

\begin{figure}[H]
  \centering
  \includegraphics[width=1.0\linewidth]{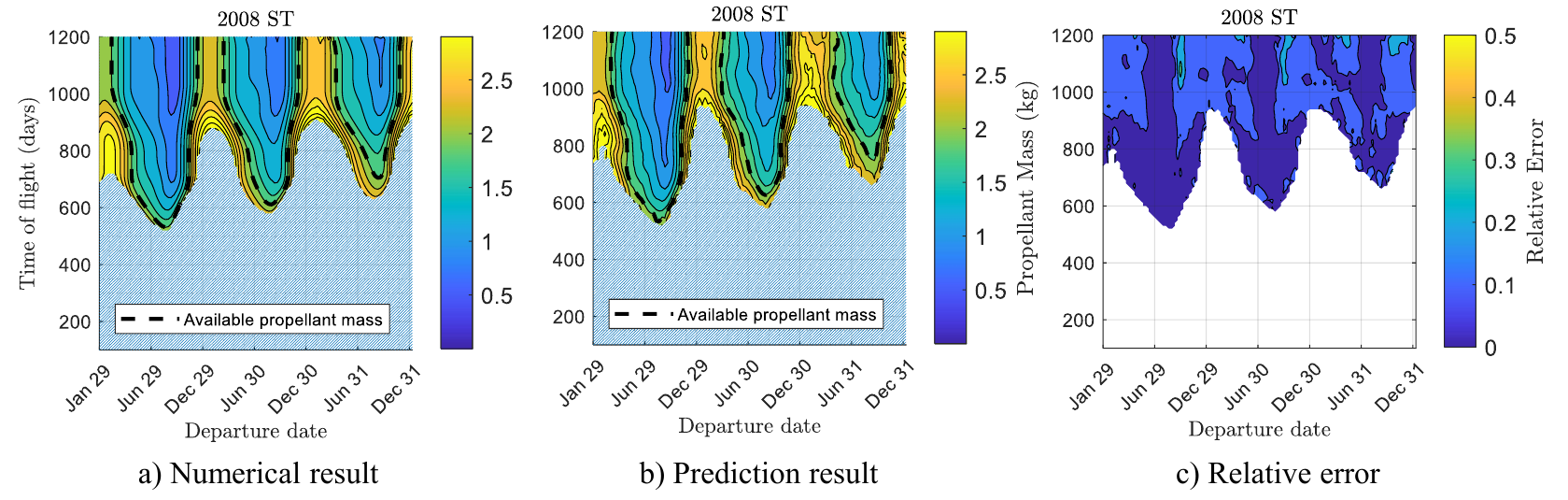}
  \caption{Porkchop plot for the Asteroid 2008 ST without variable $I_{\rm sp}$, $T_{\max}$, and path constraints.}
  \label{fig:2008ST_noconstraint}
\end{figure}

\begin{figure}[H]
  \centering
  \includegraphics[width=1.0\linewidth]{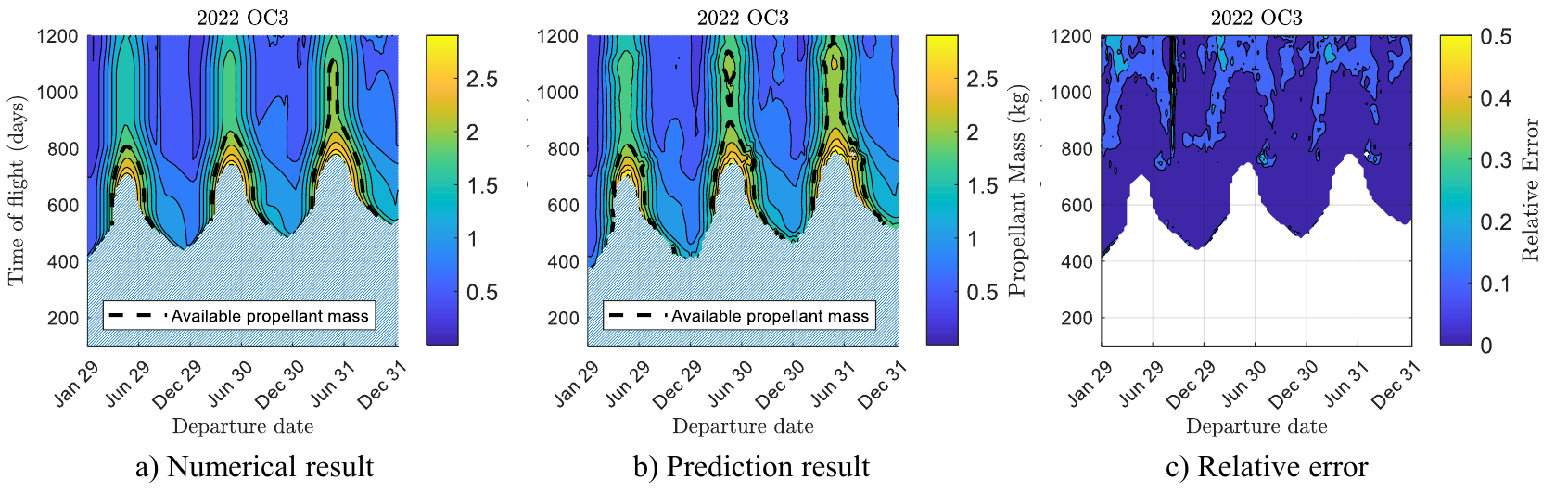}
  \caption{Porkchop plot for the Asteroid 2022 OC3 without variable $I_{\rm sp}$, $T_{\max}$, and path constraints.}
  \label{fig:2022OC3_noconstraint}
\end{figure}

This section further investigates the causes of the non-negligible errors
observed in the previous section. Intuitively, these errors may arise from the
approximation inherent in the neural network model or from the complexity of
the orbital dynamics and constraint. To isolate the effect of the neural
network, variable specific impulse, thrust, and path constraints are removed in
this experiment, leaving only the core dynamical model. The resulting relative
error plots are shown in Figures~\ref{fig:2012LA_noconstraint},
~\ref{fig:2008ST_noconstraint}, Fig.~\ref{fig:2022OC3_noconstraint}. As
illustrated in the figures, the neural network consistently captures the
overall trends and salient features of the trajectory across all asteroid
scenarios. This demonstrates the model’s capability to represent the underlying
orbital dynamics and validates the effectiveness of the proposed learning-based
approach. The discrepancies between the porkchop plots and the neural network
predictions for asteroid 2008 ST are further analyzed to investigate the
underlying reasons for the observed differences. As shown in
Figure~\ref{fig:2008ST_noconstraint}, removing the constraints leads to
significantly different porkchop plot structures compared to the original ones
in Figure~\ref{fig:porkchop_2008ST}. This indicates that the new structures
arise from the inclusion of variable specific impulse $I_{\rm sp}$, maximum
thrust $T_{\max}$, and path constraints specific to this scenario.

\begin{figure}[H]
  \centering
  \includegraphics[width=0.8\linewidth]{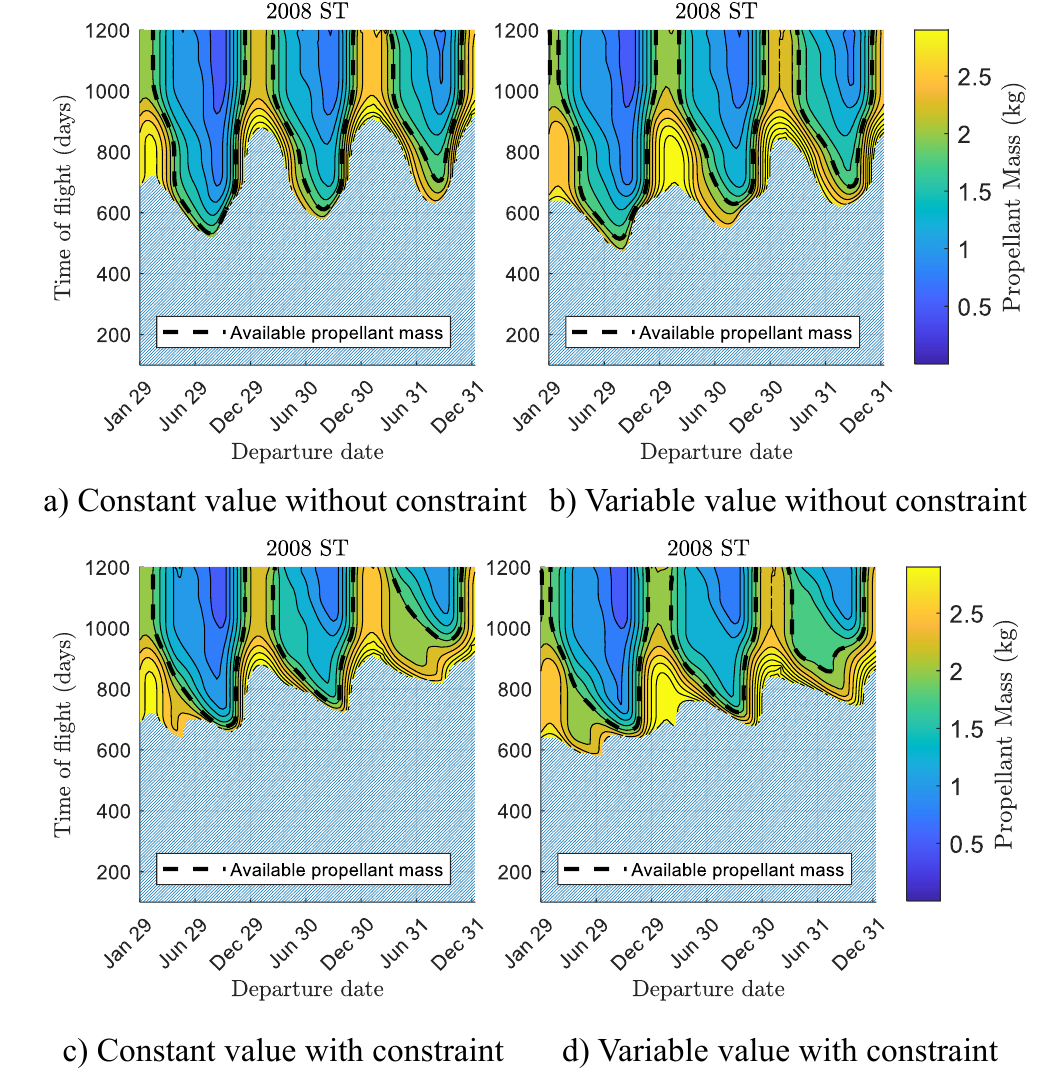}
  \caption{Comparison of porkchop plots with variable $I_{\rm sp}$, $T_{\max}$ and path constraints by optimal control method.}
  \label{fig:compare_constraints}
\end{figure}

Therefore, the individual effects of variable specific impulse, maximum thrust,
and path constraints on the porkchop plots are further investigated. To this
end, porkchop plots with and without variable specific impulse, maximum thrust,
and path constraints are generated using a direct method of optimal control, as
shown in Fig.~\ref{fig:compare_constraints}. The variations in specific impulse
and maximum thrust have a relatively minor impact on the orbital transfer,
resulting in only slight structural changes in the porkchop plots. In contrast,
the introduction of path constraints significantly affects the porkchop plot,
particularly near its boundaries, and substantially alters the overall plot
structure. This suggests that, in the later stages of mission planning, the
influence of path constraints must be carefully accounted for to optimize
launch window selection.

\begin{figure}[H]
  \centering
  \includegraphics[width=1.0\linewidth]{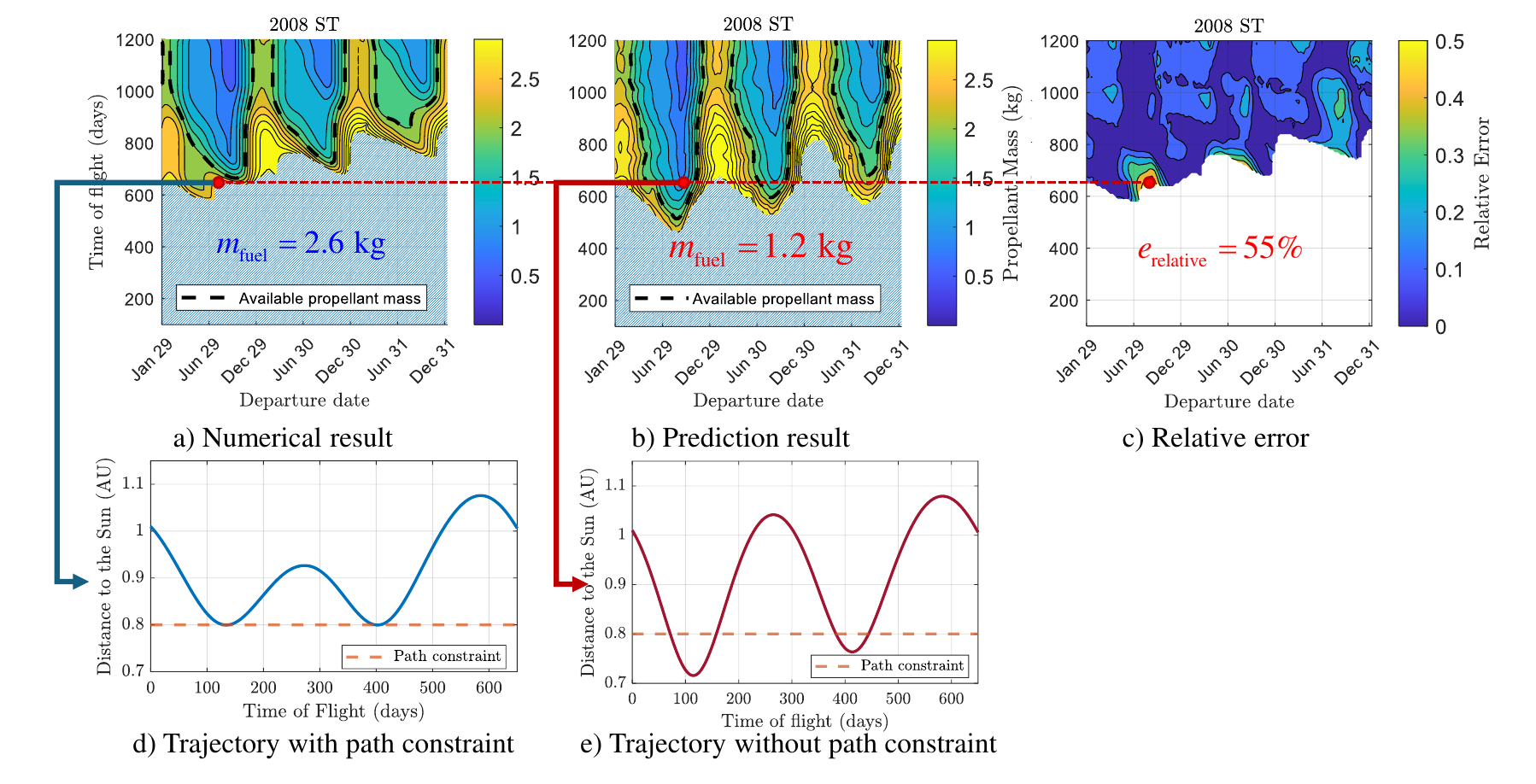}
  \caption{Porkchop plot for the 2008 ST Asteroid with path constraint.}
  \label{fig:porkchop_2008ST_with_path}
\end{figure}

Further analysis of the reverse-computed trajectories, shown in Fig.
\ref{fig:porkchop_2008ST_with_path}, reveals that the optimal control method,
constrained by path boundaries, produces trajectories that closely follow the
0.8 AU limit. In contrast, the neural network approach, which does not impose
path constraints, results in infeasible trajectories in certain regions. This
shows that at the point with the largest prediction error, there is a
significant difference between subplots (d) and (e), which correspond to the
path-constrained scenarios. This indicates that the current neural network
model may produce large estimation errors when physical path constraints are
present. This also highlights the need for future research to incorporate such
constraints into the network model.

\begin{figure}[H]
  \centering
  \includegraphics[width=0.8\linewidth]{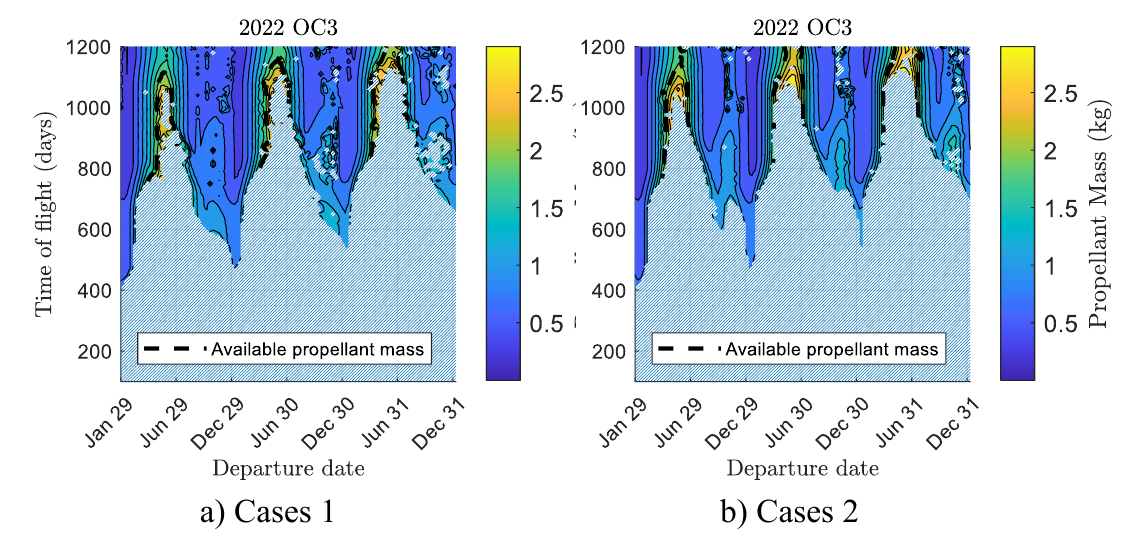}
  \caption{Example of failed cases of solution continuation with direct method, due to discontinuities in $I_{\rm sp}$, $T_{\max}$.}
  \label{fig:porkchop_2022OC3_failed1}
\end{figure}

In some cases, the presence of multiple local optima causes the optimal control
method based on the solution-continuation technique to fail in capturing the
complete global cost landscape, as illustrated in
Fig.~\ref{fig:porkchop_2022OC3_failed1}. Although this issue can be mitigated
through extensive trial runs, such an approach would inevitably lead to
significantly increased computation time. This challenge becomes particularly
pronounced when screening a large number of candidate asteroids for
exploration, where the global computational burden becomes increasingly
prohibitive. This scenario presents a compelling case for the use of neural
network models. When trained on a sufficient amount of high-quality data,
neural networks can leverage their global approximation capability to provide
smooth and reasonable fuel consumption estimates. Moreover, due to their high
computational efficiency, neural networks enable the use of global optimization
algorithms within a reasonable time frame, even in scenarios requiring further
optimization like tuning the launcher's velocity in this paper's study, making
them a promising alternative.

\section{Conclusion}

This work presents a comparative study of the performance of optimal control
methods and neural network-based estimation approaches on porkchop plots, a
widely adopted tool in mission analysis. In simple cases without path
constraints, the neural network approach is able to reproduce the structural
characteristics of the original porkchop plot while maintaining relative
prediction errors below 10\% in the majority of regions. In scenarios involving
multiple local optima where continuation-based optimal control methods may
fail, the neural network approach remains stable due to its independence from
initial guesses. This facilitates efficient identification of asteroid targets
with scientific or engineering value by mission designers. Currently, practical
engineering constraints, such as path constraints, have not yet been
incorporated into the neural network model, which limits its applicability in
some cases. Future research could focus on integrating these constraints into
the learning framework, or exploring approaches guided by physical constraints.
This paper also presents a new perspective on the application of neural
networks in astrodynamics. The generation of porkchop plots is inherently
intuitive and easy to visualize, while enabling the validation and verification
of a large number of orbital transfer scenarios, thus allowing for a
straightforward assessment of prediction quality. More importantly, this tool
is widely used in mission preliminary design and is particularly valuable in
early-stage planning, where a higher tolerance for solution inaccuracy is often
acceptable. The introduction of neural networks can significantly enhance the
efficiency of mission designers. Furthermore, due to their independence from
domain-specific expertise, planetary science teams can quickly perform
preliminary filtering during the initial mission phase, thereby reducing
communication and iteration costs, which offers practical value for engineering
applications.

\vspace{6pt}

\authorcontributions{Conceptualization, Z.Z.; methodology, Z.Z.; software, Z.Z., N.M. and G.O.P.; validation,  Z.Z., N.M., G.O.P. and Y.Z.; writing---original draft preparation, Z.Z.; supervision, F.T.; project administration, F.T.; All authors have read and agreed to the published version of the manuscript.}

\funding{This research was funded by the National Natural Science Foundation of China grant number 12532018.} 





\conflictsofinterest{The authors declare no conflict of interest.}



\reftitle{References}
\externalbibliography{yes}
\bibliography{references}

\end{document}